\documentclass{article}
\usepackage{amssymb, verbatim}
\usepackage{hyperref}
\usepackage{amsmath}

\def\Ric{\mathop{\rm Ric}}

\def\dist{\mathop{\rm dist}}

\def\Riem{\mathop{\rm Rm}}

\def\supp{\mathop{\rm supp}}

\def\RR{\mathop{\mathbb{R}}}
\def\HH{\mathop{\mathbb{H}}}
\def\EE{\mathop{\mathbb{E}}}

\def\CC{\mathop{\mathbb{C}}}

\def\GG{\mathop{\mathcal{G}}}

\def\Cap{\mathop{\rm Cap}}

\def\be{\begin{eqnarray}}
\def\ee{\end{eqnarray}}
\def\beg{\begin{eqnarray*}}
\def\ees{\end{eqnarray*}}

\newcommand{\qed}{\hfill$\Box$}
\newtheorem{theorem}{Theorem}[section]
\newtheorem{proposition}[theorem]{Proposition}
\newtheorem{lemma}[theorem]{Lemma}

\newtheorem{corollary}[theorem]{Corollary}

\date{\small\it October 5, 2018}
\author{Brian Weber}
\title{Topology of K\"ahler manifolds with weakly pseudoconvex boundary}

\begin{document}

\maketitle

\begin{abstract}
We study K\"ahler manifolds-with-boundary, not necessarily compact, with weakly pseudoconvex boundary, each component of which is compact.
If such a manifold $K$ has $l\ge2$ boundary components (possibly $l=\infty$), then it has first betti number at least $l-1$, and the Levi form of any boundary component is zero.
If $K$ has $l\ge1$ pseudoconvex boundary components and at least one non-parabolic end, the first betti number of $K$ is at least $l$.
In either case, any boundary component has non-vanishing first betti number.
If $K$ has one pseudoconvex boundary component with vanishing first betti number, the first betti number of $K$ is also zero.
Especially significant are applications to K\"ahler ALE manifolds, and to K\"ahler 4-manifolds.
This significantly extends prior results in this direction (eg. Kohn-Rossi), and uses substantially simpler methods.
\end{abstract}

\section{Introduction}

On K\"ahler manifolds with weakly pseudoconvex boundary, we use certain connections between pseudoconvexity and harmonic function theory to obtain topological constraints on both the manifold and its boundary.
For the purposes of our work, we take a boundary component $L\subset\partial{K}^m$ of a complex manifold $(K^m,J)$ to be weakly pseudoconvex if it has a plurisubharmonic defining function, meaning a differentiable defining function $f\le0$ with $\sqrt{-1}\partial\bar\partial\,f\ge0$, and is strongly pseudoconvex if $\sqrt{-1}\partial\bar\partial{f}>0$.
We consider aspects of both the real and complex geometry of $K^m$, so it will be convenient to sometimes use $-\frac12dJd\varphi=\sqrt{-1}\partial\bar\partial{}f$.

A Green's function on a complete manifold is any function $G$ defined on the complement of a point $x$ with $\triangle_g{G}=-\delta_x$ distributionally; here $\triangle_g=g^{ij}\nabla^2_{\partial_i,\partial_j}$ is the rough Laplacian.
A complete manifold is called non-parabolic if it admits a Green's function that is bounded on one side, and parabolic if not.
This definition applies to manifolds-with-boundary, assuming the boundary is compact, by requiring Neumann boundary conditions; therefore ends of manifolds (connected unbounded domains with compact boundary) may themselves be referred to as parabolic or non-parabolic.
It is known that a complete manifold with one non-parabolic end is non-parabolic.

The K\"ahler metric on a neighborhood of a pseudoconvex boundary component can be extended to make a complete K\"ahler end, by choosing an appropriate potential function.
We show that the metric on this end is non-parabolic in a strong sense; in the terminology of Section \ref{SectionParaNonPara}, an end formed this way is {\it distinguishable}.
Thus the theory of non-parabolic ends can be applied near any compact pseudoconvex boundary component of any K\"ahler manifold.

Non-constant but uniformly bounded harmonic functions exist on Riemannian manifolds when they have two or more non-parabolic ends (eg. \cite{LT1}), so after extending the pseudoconvex boundaries of $K$ as described, we obtain bounded non-constant harmonic functions.
A simple argument shows that an harmonic function $h$ obtained in this way is actually pluriharmonic: $\partial\bar\partial{h}=0$.
A particular consequence is that all boundary components have a defining function with zero Levi form.
In addition, $\partial{h}$ clearly carries non-trivial Dolbeault cohomology in $H^{1,0}(K)$.
Less trivially, we also show that $Jdh$ carries non-trivial de Rham cohomology.

Throughout, the hypotheses on our manifolds are the following:
\begin{itemize}
\item[$(*)$] \hfill\parbox{.95\linewidth}{$(K^m,\,J,\,\omega_0)$  is a K\"ahler manifold-with-boundary of complex dimension $m$, with $n$ many non-parabolic ends (possibly $n=0$ or $\infty$).
If $\{L_i\}_{i=1}^l$ are its boundary components (possibly $l=\infty$), then each $L_i$ is compact, smooth, and has a defining function $f_i$ defined in a neighborhood $U_i$ of $L_i$ so that $\sqrt{-1}\partial\bar\partial{f}_i\ge0$ on $U_i$.
We require the $U_i$ be disjoint, and that a constant $\epsilon>0$ exist so that each $f_i$ satisfies $|df_i|>\epsilon$, and so that $U_i$ contains an $\epsilon$-tubular neighborhood around $L_i$.}
\end{itemize}
In short, our manifolds are weakly pseudoconvex with compact boundary components, and, when there are infinitely many boundary componenets, a uniformity property on the gradients of the defining functions and and on the sizes of their domains of definition.
Our conclusions are that under the condition ($*$), the topology of $K$, the topology of the $L_i$, and the CR structure of the $L_i$ have some constraints.
Our main technical result is the following:
\begin{proposition}[{\it{cf}}. Proposition \ref{PropConstructionOfFunctions}] \label{PropPluriharmonics}
Assume $(K^m,J,\omega_0)$ satisfies $(*)$ and $l+n\ge2$.
Then for each $i\le{l}$, a non-constant pluriharmonic function $h_i:K\rightarrow[0,1]$ exists with $h_i=1$ on $L_i$ and $h_i=0$ on $L_j$ for $j\ne{i}$.
\end{proposition}
The notion of pluriharmonicity depends just on the complex structure, but noteworthy is that a pluriharmonic function is harmonic with respect to any compatible K\"ahler metric.
We use the $h_i$ from Proposition \ref{PropPluriharmonics} to prove that some betti numbers are positive.
Let $b^{p,q}(K)=\dim_{\CC}H^{p,q}_{DB}(K)$ and $b^p(K)=\dim_{\RR}H^p_{DR}(K)$ be the dimensions of the respective Dolbeault and de Rham groups.
\begin{theorem}[cf. Theorem \ref{ThmNonTrivDeRhamClass}] \label{ThmBettiNoManifold}
Assume $(K^m,J,\omega_0)$ satisfies $(*)$ and that $l\ge1$.
Then $b^{1,0}(K)\ge{l}-1$ and $b^1(K)\ge{l}-1$.
If in addition $n\ge1$, then $b^{1,0}(K)\ge{l}$ and $b^1(K)\ge{l}$.
\end{theorem}
\begin{theorem}[cf. Theorem \ref{ThmBoundaryBetti}] \label{ThmBettiNoEnds}
Assume $(K^m,J,\omega_0)$ satisfies $(*)$ and $l\ge1$.
If $l+n\ge2$, then each pseudoconvex boundary component $L_i$ has $b^1(L_i)\ge1$, and has a pluriharmonic defining function (in particular, each boundary component is Levi-flat).
\end{theorem}

In the case of complex dimension 1, Theorem \ref{ThmBettiNoEnds} is obvious and Theorem \ref{ThmBettiNoManifold} is not much more difficult.
Of course if $K^1$ is any complex 1-manifold with non-trivial boundary then $H^{1,0}(K)\ne0$, since any non-constant harmonic function $h$ provides a representative (namely $\partial{h})$  of a non-trivial $H^{1,0}$ class.
Less trivially, Theorem \ref{ThmBettiNoManifold} says $\dim{H}^1(K)\ge{l}-1$, although in dimension 1 this can be proved with a relative homology sequence.

Nevertheless it is instructive to see how the proofs of Theorems \ref{ThmBettiNoManifold} and \ref{ThmBettiNoEnds} work in dimension 1, as the general case is no more difficult once Proposition \ref{PropPluriharmonics} is accepted.
Let $K^1$ be a compact, complex $1$-manifold with smooth boundary components $\{L_i\}_{i=1}^l$ (these are {\it automatically} pseudoconvex).
Let $h_i$ be the harmonic function with $h_i=1$ on $L_i$ and $h_i=0$ on $L_j$ when $j\ne{i}$.
In the 1-dimensional case we have $2\sqrt{-1}\partial\bar\partial{f}=\triangle{f}$ for functions $f$, so harmonic functions are pluriharmonic.
Since also $-dJd{f}=2\sqrt{-1}\partial\bar\partial{f}$, each of the 1-forms $Jdh_i$ represents a class in $H^1(K)$.

To prove that this class is non-trivial, assume on the contrary that a function $f_i$ exists with $-Jdh_i=df_i$.
A computation shows that the function $z_i={h}_i+\sqrt{-1}f_i$ is holomorphic, and sends $K$ to the strip $\{0\le{Re}(z_i)\le1\}\subset\CC$.
The boundary of $K$ is mapped to the union of lines $\{Re(z_i)=0\}\cup\{Re(z_i)=1\}$, and by the open mapping theorem the image of $z_i$ has no other boundary.
However $z_i$ has no poles (as $dh_i$ and therefore $dz_i$ are bounded), so the image of $z_i$ in $\CC$ is compact, has non-empty interior, and has boundary within the parallel lines $\{Re(z_i)=0\}\cup\{Re(z_i)=1\}$.
Since this is an impossibility, we conclude one cannot solve $-Jdh_i=df_i$ for $f_i$; therefore $[Jdh_i]\in{H}^1(K)$ is a non-trivial class.
For Theorem \ref{ThmBettiNoEnds}, simply note that by restricting $Jdh_i$ to a collar neighborhood of $L_i$ and applying the same argument, we obtain a non-trivial class in $H^1([0,\epsilon]\times{L}_i)\approx{H}^1(L_i)$.

We present a few corollaries of our main theorems.
\begin{corollary} \label{CorStructureKALE}
Assume $(K^m,J,\omega_0)$ satisfies $(*)$ and has one boundary component $L$ with $\pi_1(L)$ finite.
Then $L$ is the only boundary component of $K$, all ends of $(K,J,\omega_0)$ are parabolic, and $b^1(K)=0$.
If in addition $(K^m,J,\omega_0)$ has no parabolic ends, then $b^{2m-1}(K)=0$.
If $K$ has complex dimension $2$ and no parabolic ends, then $\chi(K)\ge1$.
\end{corollary}
\underline{\sl Pf} \; If a non-parabolic end exists or if another pseudconvex boundary component exists, Theorem \ref{ThmBettiNoEnds} implies $b^1(L)>0$, contradicting the finiteness of $\pi_1(L)$.
If $b^1(K)\ne0$ then we can pass to the universal cover $\widetilde{K}$ of $K$, with the lifted K\"ahler structure $(\widetilde{J},\widetilde{\omega}_0)$.
It is easily seen that $\widetilde{K}$ continues to satisfy $(*)$.
Since $\pi_1(L)$ is finite, each component of the pre-image of $L$ is compact.
Therefore the pre-image of $L$ has infinitely many components, each of which is compact.
Letting $\widetilde{L}\subset\partial\widetilde{K}$ be one of these components, Theorem \ref{ThmBettiNoEnds} applied to $\widetilde{K}$ implies that $b^1(\widetilde{L})>0$, again an impossibility, proving that $b^1(K)=0$.

Finally assume $K$ has no parabolic ends; then $K$ has no ends.
Poincare duality gives $H^1(K)\approx{H}^{2m-1}(K,L)$ and $H^{2m-2}(L)\approx{H}^1(L)=\{0\}$, so the relative homology sequence gives $H^{2m-1}(K,L)\approx{H}^{2m-1}(K)$.
Therefore $b^{2m-1}(K)=b^{2m-1}(K,L)=b^1(K)=0$.
In the 2-dimensional case this means $b^1(K)=b^3(K)=0$ so
\begin{equation}
	\chi(K)\;=\;1\,-\,b^1(K)\,+\,b^2(K)\,-\,b^3(K)
\;=\;1\,+\,b^2(K) \;\ge\;1.
\end{equation}
\qed

An end of a Riemannian manifold is called asymptotically locally Euclidean (ALE) if it is diffeomorphic to a quotient of $\RR^k\setminus{B}(1)$ by a finite subgroup of $O(k)$ (or of $U(k/2)$ in the K\"ahler case), and also has $|\Riem|=o(r^{-2})$ where $r$ is the distance to some fixed point.
Theorem \ref{ThmBettiNoEnds} can be used to show that a K\"ahler manifold of complex dimension at least 2 (whether it is of finite type or not) that has an ALE end has {\it only one} ALE end.
Thus we recover a well-known result implied by the statement of Theorem 4.2 of \cite{LT1}, and, assuming $K$ has finite type, by the theorems of Kohn-Rossi \cite{KR} and Kohn \cite{K1}.
These past results are discussed in the remarks below.
\begin{corollary} \label{CorKALE}
Assume $(K,J,\omega_0)$ has complex dimension at least 2, satisfies $(*)$, and has an ALE end.
Then every other end of $K$ is parabolic, and $H^1(K)=0$.
If $K$ has no parabolic ends, then $H^{n-1}(K)=0$.
If $K$ has no parabolic ends and complex dimension 2, then $\chi(K)\ge1$.
\end{corollary}
\underline{\sl Pf} \; We can assume $K'\subset{K}$ is an ALE end so that the boundary of $K\setminus{K}'$ is diffeomorphic to a quotient of an $(n-1)$-sphere and is geometrically locally convex, and therefore pseudoconvex.
Corollary \ref{CorStructureKALE} applied to $K\setminus{K}'$ then provides the conclusion.
\qed

{\bf Remark.} Our results substantially expand what was previously known.
In particular, we require just non-negativity of eigenvalues of the Levi form instead of positivity, and we do not require that the manifold have compact closure.
In trade, we require the manifold be K\"ahler rather than just Hermitian.
Our new hypotheses, particularly our allowance of infinitely many boundary components, allow us to pass to universal covers and still use our main theorems.
This allows us to prove our statements on the betti numbers of the boundary components themselves, as well as the statements about the Euler characteristic in Corollaries \ref{CorStructureKALE} and \ref{CorKALE}.

{\bf Remark.} The effect of boundary pseudoconvexity on cohomology has been studied extensively by many authors.
Hilbert space methods were developed in Kohn \cite{K1} \cite{K2}, Andreotti-Vesentini \cite{AV}, and H\"ormander \cite{Ho1} for the purpose of solving $\bar\partial$-Neumann problems and non-homogeneous $\bar\partial$-problems.
One result was a proof that the pseudoconvexity of subdomains of $\CC^n$ or of Stein manifolds gives rise to strong cohomological vanishing theorems.
In addition, a Hodge decomposition on compact complex manifolds-with-boundary holds in a given bidegree provided the boundary satisfies a certain pseudoconvexity condition (which in any bidegree is implied by strong pseudoconvexity; see \cite{KR}).
That is, if $\triangle$ is the $\bar\partial$-Laplacian, there is a compact operator $G$ so that 
\begin{eqnarray}
{\bigwedge}^{p,q}=\triangle\circ{G}\left({\bigwedge}^{p,q}\right)\oplus\,\mathcal{H}^{p,q}, \label{EqnStrongHodge}
\end{eqnarray}
where $\bigwedge^{p,q}$ is the space of $C^\infty$ forms of the indicated bidegree, and $\mathcal{H}^{p,q}$ is the space of harmonic $(p,q)$-forms.
The delicate analysis required for the proofs involves H\"ormander-style estimates with plurisubharmonic functions.

{\bf Remark.} Kohn-Rossi \cite{KR} used the Hodge decomposition (\ref{EqnStrongHodge}) to solve a number of boundary value problems, one of which was the following:
if $K$ is a compact complex manifold whose boundary satisfies an appropriate convexity condition (the Levi form is positive for instance; for the precise condition see 7.1 of \cite{KR}), if $f$ is a function on $\partial{K}$ that satisfies a certain compatibility condition (namely that $\bar\partial_bf=0$, where $\bar\partial_b$ is the restriction of the $\bar\partial$-operator to the boundary), and if $f$ is orthogonal to the restriction of $\mathcal{H}^{m,m-1}$ to the boundary, then $f$ is the restriction to $\partial{K}$ of a holomorphic function $F$ on $K$.
A corollary, described in the next paragraph, is that a {\it compact} Hermitian manifold, all of whose boundary components are strictly pseudoconvex, must have connected boundary.
Assuming the manifold is K\"ahler, this is a weakened version of our Theorem \ref{ThmBettiNoEnds}.

To see how Kohn-Rossi proved this, first note that strict pseudoconvexity implies $\mathcal{H}^{m,m-1}$ is finite dimensional, by the previous remark.
The existence of a single non-constant holomorphic function on $A$ is implied by Theorem 9.1 of \cite{K1} (see also \cite{Gr}).
By taking powers of this function, we see that the vector space of holomorphic functions is infinite dimensional.
Consider the subspace of holomorphic functions $span_{\mathbb{C}}\{A^i\}_{i=0}^\infty$ spanned by powers $A$.
To any function $F=\sum{}c_iA^i$ in this subspace, construct the the function $g:\partial{M}\rightarrow\CC$ by multiplying $F|_{\partial{M}}$ by different constants on each component of $\partial{M}$.
By the finite-dimensionality of $\mathcal{H}^{m,m-1}$, we can choose the $c_i$ so that $g$ is orthogonal to the restriction of $\mathcal{H}^{m,m-1}$ to $\partial{M}$.
Thus, since also $\bar\partial_bg=0$, we can extend $g$ to a holomorphic function $G$ on $K$, by the Kohn-Rossi extension theorem mentioned above.
Then $G/F$ will be a meromorphic function that is locally constant on the boundary, and therefore constant by unique continuation.
But since $G/F$ takes different values on each boundary component, the boundary must have just one component.

We have presented this summary of a key part of the Kohn-Rossi method (adapted from 7.2 and 7.3 of \cite{KR}) in order to emphasize its reliance on the finite dimensionality of $\mathcal{H}^{m,m-1}$, which is not required in our setting.

{\bf Organization.} In section \ref{SectionParaNonPara} we collect the material on harmonic function theory used in the proof of Proposition \ref{PropPluriharmonics}.
Our main concern is with exactly how harmonic functions with 2-sided bounds are constructed---the method is termed compact exhaustion; see \cite{LT0} and \cite{LT1}.
We also introduce the useful notion of {\it distinguishability}, which is a strengthened form of non-parabolicity that appears in our context, and we conclude with an example showing that distinguishability is strictly stronger than non-parabolicity.
Section \ref{SectionFunctionKahler} contains the proofs of Proposition \ref{PropPluriharmonics} and of Theorems \ref{ThmBettiNoManifold} and \ref{ThmBettiNoEnds}.

{\bf Acknowledgements}. The author would like to thank Xiuxiong Chen and Claude LeBrun for several useful conversations, and Charles Epstein for making him aware of the results of Kohn-Rossi \cite{KR}.
Special thanks go to New York University's Courant Institute, which provided working space for the author during the writing of this paper, and to the National Science Foundation\footnote{Any opinions, findings and conclusions or recommendations in this material are those of the author(s) and do not necessarily reflect the view of the National Science Foundation} for the grant DMS-0635607002 that provided support.

\section{Parabolic and non-parabolic ends of Riemannian manifolds} \label{SectionParaNonPara}

The literature on harmonic function theory on complete manifolds is very large.
Here we gather some well-known results that will be useful later, and introduce the notion of the {\it distinguishability} of an end, which means, roughly speaking, that the end can be separated from the rest of the manifold by a bounded harmonic function.
We show that this notion is strictly stronger than non-parabolicity.
In this section we are concerned only with Riemannian, not K\"ahler, structures.

A function $\GG_x:M\rightarrow\RR$ is called a Green's function at the point $x$ if $\triangle\GG_x=-\delta_x$ in the sense of distributions.
A complete manifold is called parabolic if it admits no positive Green's function, and non-parabolic otherwise.
These definitions are equally good on manifolds with compact boundary, with Green's functions made to satisfy Neumann conditions on boundary components.
With $c_{n-1}$ the area of the unit $(n-1)$-sphere, the Green's functions (at the origin) of the flat manifolds $\RR^n$ are
\begin{eqnarray}
{\GG}_x(y) &=& \begin{cases}
-\frac{1}{2\pi}\log|x-y| & if \; n=2    \\
\frac{1}{(n-2)c_{n-1}}|x\,-\,y|^{2-n}  & if \; n\ne2.
\end{cases} \label{EqnFlatGreens}
\end{eqnarray}
Therefore $\RR^n$ is parabolic when $n=2$ and non-parabolic when $n>2$.

Related to parabolicity is the notion of capacity.
Given any set $\Omega\subset{M}$ with compact closure, we define its capacity $\Cap(\Omega)$ to be an infimum of Dirichlet integrals:
\begin{eqnarray}
\Cap(\Omega) &=&\inf_\varphi \int_M|\nabla\varphi|^2 \label{EqnDefCap}
\end{eqnarray}
where the infimum is over all $\varphi\in{C}_c^{0,1}(M)$ with $\varphi\ge1$ on $\Omega$.
Assuming $\Omega$ is a smooth domain and $\Cap(\Omega)>0$, the infimum is obtained by a Lipschitz function $\varphi$ with $\varphi=1$ on $\Omega$, $\triangle\varphi=0$ outside $\Omega$, and $\varphi\rightarrow0$ along some (but not necessarily every) sequence of points that diverges to infinity.
If $\Cap(\Omega)=0$, a minimizing sequence will converge to a constant function.
The connection between capacity and parabolicity is the following proposition, which can be found for instance in \cite{H1}, and also follows from (2) of Proposition 1.2 in \cite{LT2}.
\begin{proposition}
A Riemannian manifold $(M,g)$ is non-parabolic if and only if it has a sub-domain with compact closure and positive capacity.
\end{proposition}
\qed

A geometric phenomenon totally absent on $\RR^n$, $n\ne1$, is the possibility of separating unbounded sets with domains of compact closure; this gives rise to the notion of ends.
If $\Omega$ is a pre-compact domain, we call any unbounded component of $M\setminus\overline{\Omega}$ an {\it end} of $M$ with respect to $\Omega$.
We shall call any connected, unbounded subset $M'$ an end if $\partial{M'}$ is non-empty and compact---we usually leave the domain $\Omega$ implicit.
Note that an end may have two or more non-intersecting subsets that are themselves distinct ends.

Capacity, and therefore parabolicity and non-parabolicity, is a notion that can be attributed to an end, not just to a manifold.
If $M'$ is an end, we define its capacity to be
\begin{eqnarray}
\Cap(M') &=& \inf_\varphi\int_{M'}|\nabla\varphi|^2 \label{EqnDefEndCap}
\end{eqnarray}
where the infimum is taken over all $C^{0,1}(M')$ functions $\varphi$ of compact support with $\varphi=1$ on $\partial{M'}$.
If $\partial{M}'$ is smooth and $\Cap(M')>0$, the capacity is realized by some harmonic $C_0^\infty(M')$ function $\varphi$ with $\varphi=1$ on the boundary.
We call an end $M'$ non-parabolic if $\Cap(M')>0$, and parabolic if $\Cap(M')=0$.
Clearly a manifold is non-parabolic if it has one non-parabolic end.
There is also the following alternative characterization ({\it cf}. Proposition 1.2 of \cite{LT2}).
\begin{lemma}
 \label{LemmaNonParaBarrier}
An end $M'$ of $M$ is non-parabolic if and only if a superharmonic function $f$ exists on $M'$ with $\inf_{\partial{M'}}f>0$ and $f\rightarrow0$ along some sequence of points that diverges to infinity.
\end{lemma}
\underline{\sl Pf} \; The proof here is similar to that in \cite{LT2}; we go through it because some details will be used later.
First assume the stated function $f$ exists.
After multiplying by $(\inf_{\partial{M}'}f)^{-1}$ we can assume $f\ge1$ on $\partial{M}'$.
Let $\varphi_i$ be a minimizing sequence for (\ref{EqnDefEndCap}); then $\varphi_i\ge1$ on $\partial{}M'$ and the support $M'_i\subset{}M'$ of $\varphi_i$ is compact---indeed we can assume $\{M'_i\}$ is an exhaustion of $M'$.
A simple argument (which we omit) states that we can replace $\varphi_i$ by $\min\{1,\varphi_i\}$ to obtain a function with smaller Dirichlet integral, and that we can replace $\varphi_i$ with an harmonic function with the same boundary values, and also obtain a function with strictly smaller Dirichlet integral (since harnomic functions minimize the Dirichlet energy among all functions with given boundary data).
Therefore we may take the $\varphi_i$ to be harmonic, have compact support, satisfy $\varphi=1$ on $\partial{M}'$, and that $M'_i=\supp\varphi_i$ is an exhaustion of $M'$ by compact sets.
Since $f$ is superharmonic and $f\ge\varphi_i$ on each $\partial{}M'_i$, we have $\varphi_i\le{f}$.
Since $0\le\varphi_i\le1$, a subsequence will converge to an harmonic function $\varphi$, and we retain $f\ge\varphi$.
The Dirichlet integrals $\int|\nabla\varphi_i|^2$ decrease monotonically and converge to $\int|\nabla\varphi|^2$, so that $\varphi$ is non-constant and has a finite (but non-zero) Dirichlet integral.

For the converse, now assume $M'$ is non-parabolic.
We may assume $M'$ has a smooth boundary, as shrinking $M'$ increases its capacity.
Letting $M'_i$ be a compact exhaustion of $M'$ so that $\partial{M}'\subset\partial\Omega_i$, let $\varphi_i$ be harmonic functions with $\varphi_i=1$ on $\partial{M}'$ and $\varphi_i=0$ on $\partial{}M'_i\setminus\partial{M}'$.
We have $\int_{\partial{M}'}\frac{\partial\varphi_i}{\partial\hat{n}}=\int_{M'}|\nabla\varphi_i|^2>\Cap(M')>0$, and (due to the Hopf lemma) that $\int|\nabla\varphi_i|^2$ decreases monotonically.
The $\varphi_i$ are bounded and harmonic, so certainly they converge along a subsequence; indeed under our construction the $\varphi_i$ are monotonic with respect to $i$, so the $\varphi_i$ actually converge.
Set $\varphi=\lim_i\varphi_i$.
Because $\int_{M'}|\nabla\varphi_i|^2\ge\Cap(M')$ (by the definition of capacity), and bacause the Dirichlet integral is decreasing, Fatou's Lemma gives $\int_{M'}|\nabla\varphi|^2\ge\Cap(M')$.
Thus $\varphi$ is non-constant, and of course it is harmonic and bounded between $0$ and $1$.
By these properties and because $\varphi=1$ on $\partial{M'}$, $\varphi$ obtains a strict minimum at infinity.
We wish to prove $\inf_{M'}\varphi=0$.
Setting $\epsilon=\inf\varphi$, then $\tilde\varphi=\frac{\varphi-\epsilon}{1-\epsilon}$ is an harmonic function equal to $1$ on the boundary, and $\tilde\varphi\rightarrow0$ along some subsequence that diverges to $\infty$.
Using $\tilde\varphi$ as a barrier and following the argument of the previous paragraph, we have that actually $\varphi=\lim_i\varphi_i\le\tilde\varphi$, which means $\inf\varphi=0$.
\qed

We shall call an end $M'$ {\it distinguishable} if a positive harmonic function $\varphi$ exists on $M'$ with $\varphi=1$ on $\partial{M}'$ and $\varphi\rightarrow0$ along {\it every} sequence of points in $M'$ that diverges to $\infty$.
The following lemma is essentially obvious.
\begin{lemma} \label{LemmaDistinguishedBarrier}
An end $M'$ of a manifold is distinguishable if and only if there is a positive superharmonic function $f:M'\rightarrow\RR$ with $\inf_{\partial{M'}}f>0$ and so that $f\rightarrow0$ along every sequence of points in $M'$ that diverges to infinity.
\end{lemma}
\underline{\sl Pf} \; This follows after the constructing a harmonic function $\varphi$ as in Lemma \ref{LemmaNonParaBarrier}, by noting that (after possibly multiplying $f$ by a constant to make $f\ge\varphi$ on $\partial{M}'$), we have $f\ge\varphi>0$.
\qed

The importance of distinguishability comes from the following lemma, which states that, on a manifold-with-boundary with compact boundary, distinguishable ends can be separated from the rest of the manifold with harmonic functions, provided at least one other non-parabolic end exists.
The first assertion in the following proposition is well known (eg. \cite{LT2}).
The second assertion is new.
\begin{proposition}[Separation of distinguishable ends] \label{PropConstructionOfFunctions}
Assume $(M,g)$ is a smooth Riemannian manifold-with-boundary, with smooth boundary.
If $M$ has at least two non-parabolic ends, then there exists a non-constant harmonic function $\varphi:M\rightarrow\RR$ with $0<\varphi<1$.
If, in addition, $M'$ is a distinguishable non-parabolic end, a number $\delta'>0$ can be chosen so that if $\delta\in(0,\delta')$ and $\Omega_\delta\triangleq\{\varphi>1-\delta\}$, we have that $\Omega_\delta\subset{M}'$ and that $M'\setminus\Omega_\delta$ has compact closure.
\end{proposition}
\underline{\sl Pf} \; The method for proving the first assertion is compact exhaustion.
Namely let $f'$, $f''$ be the superharmonic barrier functions on the non-parabolic ends $M'$, $M''$, respectively, that are guaranteed by Lemma \ref{LemmaNonParaBarrier} or \ref{LemmaDistinguishedBarrier}.
Let $M_i$ be an exhaustion of $M$ by smooth, pre-compact domains, each of which separates $M'$ and $M''$.
If $\partial{M}$ is non-empty, assume $\partial{M}\subset{M}_i$.
Let $\varphi_i:M_i\rightarrow\RR$ be the harmonic function with $\varphi_i=1$ on $(\partial{M}_i\cap{M}')\setminus\partial{M}$, $\varphi_i=0$ on $\partial{M}_i\setminus{M}'$, and $\varphi_i$ satisfies von Neumann conditions on $\partial{M}$.
Set $\varphi=\lim_i\varphi_i$.
Since $\varphi_i\ge1-f'$ on $M'$ and $\varphi_i\le{f}''$ on $M''$, the same holds for $\varphi$; therefore $\varphi$ is not constant.

We can prove that on $M\setminus{M}'$, $\varphi$ is bounded strictly below $1$, unless possibly $M'$ is the only non-parabolic end.
By the maximum principle we have $\sup_{M\setminus{M}'}\varphi_i=\sup_{\partial{M'}}\varphi_i$.
Taking $i\rightarrow\infty$, the same holds for $\varphi$.
On the other hand, $\varphi\le1$ so the strong maximum principle implies that either $\varphi\equiv1$ or else $\varphi<1$ on $K$.
By the compactness of $\partial{M}'$, we have either $\varphi\equiv1$ on $K$ or else $\sup_{M\setminus{M}'}\varphi=\sup_{\partial{M}'}<1-\delta$ for all sufficiently small $\delta$.

Finally, assume $M'$ is distinguished.
We can assume the upper barrier $f'$ satisfies $f'(x_i)\searrow0$ along any sequence $x_i$ in $M'$ that diverges to infinity.
Now choose the $M_i$ so that $f'<2^{-i}$ on $M'\cap{M}_i$; clearly $M\setminus{M}_i$ is compact.
Since $\varphi\ge1-f'>1-2^{-i}$ on $M'\cap{M}_i$, we can choose a large enough $i$ so that $1-2^{-i}>\sup_{\partial{M}'}\varphi$.
Putting $\delta=2^{-i}$, the set $\Omega_\delta=\{\varphi>1-\delta\}$ is therefore a subset of ${M}'$.
Clearly $M'\setminus\Omega_\delta$ is compact as $M_i\subset\Omega_\delta\subset{M}_j$ (strict inclusion) when $2^{-j}<\delta<2^{-i}$.
\qed

We close this section with an example of an end that is non-parabolic but not distinguishable.
Let $\EE^2$ be $\RR^2$ a flat metric $g_F$ and let $\HH^2$ be $\RR^2$ with a hyperbolic metric $g_H$.
Then $\EE^2$ is parabolic with Green's function given by (\ref{EqnFlatGreens}) and $\HH^2$ is non-parabolic with Green's function $\GG_x(y)=\frac{1}{2\pi}\log\left(\frac{e^r-1}{e^r+1}\right)$, where $r=\dist(x,y)$.
Let $x_i$ (resp. $y_i$) be a sequence of points in $\EE^2$ (resp. $\HH^2$) with $x_i\rightarrow\infty$ (resp. $y_i\rightarrow\infty$), and attach $\EE^2$ to $\HH^2$ by removing small balls $B_{x_i}(\delta_i/2)$ from $\EE^2$ and $B_{y_i}(\delta_i/2)$ from $\HH^2$ ($\delta_i$ is a sequence of positive numbers), and gluing the ends of a cylinder to each pair of corresponding boundary components.
Label this manifold $(M,g)$, where $g$ is chosen so the metrics on $\EE^2\setminus\bigcup_iB_{x_i}(\delta_i)$ and $\HH^2\setminus\bigcup_iB_{y_i}(\delta_i)$ are unchanged.

The resulting manifold $(M,g)$ clearly has a single end.
That $M$ is non-parabolic follows, for instance, from Theorem 2.1 of \cite{HK} with $p=2$, after noting that the metric on the hyperbolic part of $M$ makes the volume growth of balls exponential.

Let $B_E$ be the unit ball about the origin on the Euclidean part of $M$ (we assume this does not intersect any of the $B_{x_i}(\delta_i)$), and let $M'=M\setminus{B}_E$ be the end with respect to $B_E$.
We will construct a lower barrier functions $F$ with the property that any positive harmonic function $\varphi$ on $M$ with $\varphi\ge1$ on $\partial{}B_E$ has $\varphi\ge{F}$, and that $F$ is asymptotically nonzero along some diverging sequences.

Consider the following family of functions defined a.e. on $\EE^2$:
\begin{eqnarray}
F_{\eta}(x)&=&1\,-\,\eta\log(|x|) \,+\, \sum_{a=1}^\infty{}c_a\log(|x-x_a|).
\end{eqnarray}
If the $c_a>0$ converge to zero fast enough (say $c_a=a^{-2}$ and $x_a$ has coordinates $(a,0)$), then the sum converges as $i\rightarrow\infty$.
Note that $\triangle{F}_{\eta}$ is zero aside from a delta function of weight $-\eta$ at the origin and delta functions of positive weight $c_\alpha$ at each $x_a$.
Let $D_{\eta}=\{F_{\eta}>0\}\setminus{B}_E$.
This set is pre-compact when $\eta>\eta_0\triangleq\sum_ic_i$.
We can prove that whenever $\eta>\eta_0$, we have $\varphi\ge\frac1CF_\eta$ on $D_\eta$, where $C=1+\sum_ac_a\log(1+|x_a|)$.
To see this, first note that $F_{\eta}$ reaches its maximum on $\partial{B}_E$ (assuming of course $\eta>\eta_0$) and $F_{\eta,i}$ is bounded by $C$ there; therefore $F_{\eta}$ is bounded from above by $C$ on all of $D_{\eta}$.
If we assume
\begin{eqnarray}
\delta_a\le\exp(-C/c_a),
\end{eqnarray}
it follows that $F_{\eta}\le0$ on $B_{x_a}(\delta_a)$.
Thus $F_{\eta}\le{C}\varphi$ on $\partial{D}_\eta$, so we conlcude $F_{\eta}\le{C}\varphi$ on $D_{\eta}$.

The key point is that this construction of $F_\eta$ is compatible with the gluing operation: the gluing occurs insude the balls of of radius $\delta_a$, and since $F_\eta<0$ on $\partial{}B_{x_a}(\delta_a)$ we can extend $F_\eta$ in any way on the hyperbolic part of the manifold as long as it stays negative, and $F_\eta$ will remain a lower barrier for $\varphi$ since we have assumed $\varphi>0$.

Now sending $\eta\searrow\eta_0$, we have that $F_{\eta}$ converges pointwise to $F_{\eta_0}$, and $F_{\eta_0}\le{C}\varphi$.
But  it is easily checked that $F_{\eta_0}(x)=1-\sum_ac_a\log\left(\frac{|x-x_a|}{|x|}\right)$ is asymptotically unity along most divergent sequences, privided the $c_a$ decrease quickly enough.
Therefore it is impossible that $\varphi$ is asymptotically 0.

\section{Pluriharmonic functions on K\"ahler manifolds with pseudoconvex boundary} \label{SectionFunctionKahler}

It is known that a K\"ahler metric near a pseudoconvex boundary component can be made complete by the choice of an appropriate potential function---in fact this is a defining feature; see \cite{Ohs}.
This is important enough for us that we repeat the argument here.
Let $f$ be a positive defining function for the pseudoconvex boundary component $L\subset{K}^m$.
Specifically this means three things: {\it{a}}) $f$ is defined on some neighborhood $\Omega$ of $L$, {\it{b}}) $f=0$ on $L$ and $f>0$ on $\Omega\setminus{L}$, and {\it{c}}) $f$ is pseudoconcave, so $\sqrt{-1}\partial\bar\partial{f}=-\frac12dJdf\le0$ on $\Omega$.

If $\varphi:\RR^{+}\rightarrow\RR$ is any twice-differentiable function then
\beg
dJd\,\varphi(f) &=&\varphi''(f)df\wedge{J}df \,+\,\varphi'(f)\,dJdf
\ees
Note that $df\wedge{J}df$ is always non-positive and that since $f$ is pseudoconcave, $dJdf$ is non-negative.
Thus the form
\begin{eqnarray}
\omega &=&\omega_0 \,+\,dJd\,\varphi(f)  \label{EqnDefOfOmega}
\end{eqnarray}
is positive if $\varphi''\le0$ and $\varphi'\ge0$.
If $\varphi'$ approaches infinity sufficiently quickly as $t\rightarrow0$, the corresponding metric is complete.
One obvious choice is $\varphi(t)=\log(t)$; this gives a complete manifold with constant negative bisectional curvature at infinity.
Another possibility is $\varphi(t)=-t^{-\alpha}/\alpha$ for $\alpha>0$; in this case the bisectional curvature decays to zero like $O(r^{-2})$, where $r$ is the $\omega$-distance from a fixed point.
\begin{lemma}\label{LemmaCompleteKahlerMetric}
Assume $(K^m,J,\omega_0)$ satisfies $(*)$ of the introduction.
Then a positive smooth function $f:K^m\rightarrow\mathbb{R}$ exists so that $f$ is non-strictly plurisuperharmonic, and agrees with the positive defining functions $-f_i$ on a neighborhood of $L_i$.
Further, for $\alpha\ge0$, the $(1,1)$-form
\begin{eqnarray}
\omega &=&-\frac1\alpha dJd\left(f^{-\alpha}\right) \,+\,\omega_0 \label{EqnDefOfOmega1}
\end{eqnarray}
is a K\"ahler form whose associated metric is complete near any boundary component $L_i$ (when $\alpha=0$, we take $\omega=dJd\log\,f \,+\,\omega_0$).
\end{lemma}
\underline{\sl Pf} \; Setting $f_i=\infty$ where it was otherwise undefined, by $(*)$ there is a number $\delta$ so that
\beg
f\;=\;\inf\{\delta,f_1,\dots,f_l\}
\ees
is Lipschitz (even if $l=\infty$, due to the assumption in (*) that $|df_i|>\epsilon$ on some neighborhood of $U_\epsilon$).
When $f<\delta$ then $f$ is smooth and $dJdf\ge0$.
Now we can smooth $f$ in any way that leaves it unaffected on a neighborhood of each $L_i$ by replacing $f$ by a function $\psi(f)$.
We let $\psi$ be a smooth increasing function with $\psi(t)=t$ when $t<\delta/4$, $\psi(t)=3\delta/8$ when $t>\delta/2$, and $\psi''(t)<8/\delta$ when $-\delta<t<-\delta/2$, then $\psi(f)$ is smooth, (non-strictly) plurisuperharmonic, and agrees with each $f_i$ on some neighborhood of $L_i$, as desired.

Finally we show the resulting manifold is complete.
Choose $s,S$ so $0<s<S<\delta/2$, and let $\gamma(t)$ be a path in $\{s\le{f}\le{S}\}$ from a point in $\{f=S\}$ to a point in $\{f=s\}$.
We have
\beg
|\dot\gamma|^2 &=&\omega(\dot\gamma,J\dot\gamma) \\
&=& -(1+\alpha)f^{-2-\alpha}(df\wedge{J}df)(\gamma,J\gamma) \,+\,f^{-1-\alpha}(dJdf)(\gamma,J\gamma) \,+\,\omega_0(\gamma,J\gamma) \\
&\ge& (1+\alpha)f^{-2-\alpha}(df(\dot\gamma))^2 \;=\;\frac{4(1+\alpha)}{\alpha^2}\left(\frac{df^{-\frac{\alpha}{2}}}{dt}\right)^2
\ees
where we have used $df\wedge{J}df(X,Y)=\left<\nabla{f},X\right>\left<\nabla{f},JY\right>-\left<\nabla{f},JX\right>\left<\nabla{f},Y\right>$.
We can assume $f\circ\gamma$ is $C^1$ and decreasing, so the length of $\gamma$ is estimated from below by
\beg
\int|\dot\gamma|\,dt 
\;\ge\; \frac{2\sqrt{1+\alpha}}{\alpha}\int_{f=S}^{f=s}\frac{df^{-\frac{\alpha}{2}}}{dt} \,dt
\;=\; \frac{2\sqrt{1+\alpha}}{\alpha}\left(s^{-\frac{\alpha}{2}} \,-\,S^{-\frac{\alpha}{2}}\right).
\ees
If $\alpha=0$ the appropriate expression with logarithms is obvious.
When $\alpha\ge0$, the length of $\gamma$ therefore grows unboundedly as its terminal point approaches $\partial{K}$ at $s=0$.
\qed

A real-valued function $h$ is called pluriharmonic when $\sqrt{-1}\partial\bar\partial{h}=0$.
This depends only on the complex structure, so a pluriharmonic function is harmonic with respect to any compatible K\"ahler metric.
In this section we show that the existence of more than one pseudoconvex boundary component on $(K^m,J,\omega_0)$ allows the construction of pluriharmonic functions.
We shall be careful to observe the distinction between the original metric $g_0=\omega_0(\cdot,J\cdot)$, and a choice of a complete metric $g=\omega(\cdot,J\cdot)$ given by Lemma \ref{LemmaCompleteKahlerMetric}.

Any end that comes from a pseudoconvex boundary component is distinguishable.
To see this, note that $f_i>0$ in $U_i$ and $f_i\searrow0$ along any sequence in $U_i$ that diverges in the $g$-metric, and that since both $dJdf_i$ and $\omega$ are positive we have
\be
\triangle_g{f}_i &=&\frac{-\,dJdf_i\wedge\omega^{n-1}}{\omega^n} \;<\;0
\ee
so that $f_i$ is superharmonic.
Thus the end is distinguishable by Lemma \ref{LemmaDistinguishedBarrier}.
\begin{proposition} \label{PropNonconstPluriHarmonic}
Assume $(K,J,\omega_0)$ satisfies $(*)$ of the introduction.
There exists a positive pluriharmonic function $h_i$ on $(K,J,\omega_0)$ with $h_i=1$ on $L_i$ and $h_i=0$ on $L_j$ when $j\ne{i}$.
In particular, $dJdh_i=d^{*}Jdh_i=0$.
\end{proposition}
\underline{\sl Pf} \; For each $i$ let $V_i$ be a neighborhood of $L_i$ so that $V_i$ does not intersect any boundary component of $K$ besides $L_i$.
By Propositions \ref{LemmaCompleteKahlerMetric} and \ref{PropConstructionOfFunctions}, a harmonic function $h_i$ exists on $(K,J,\omega)$ that limits to $1$ along any divergent sequence in $V_i$ and limits to $0$ along any unbounded sequence in $V_j$ for all $j\ne{i}$.

We first prove the Dirichlet integral of $h_i$ is finite.
To see this, recall how the $h_i$ are constructed: $h_i=\lim_{R\rightarrow\infty}h_{i,R}$ where $h_{i,R}$ is the harmonic function on the large ball $B_p(R)$ with $h_{i,R}=1$ on $\partial(B_p(R)\cap{V}_i)$ and $h_{i,R}=0$ on $\partial(B_p(R)\setminus{V}_i)$.
Then
\begin{equation}
\begin{aligned}
\int_{K}|\nabla{h}_{i,R}|^2 
&=\;\int_{K\setminus{V}_i}|\nabla{h}_{i,R}|^2 \,+\,\int_{{V}_i}|\nabla(1-{h}_{i,R})|^2 \\
&=\;\int_{\partial{V}_i}h_{i,R}\frac{\partial{h}_{i,R}}{\partial\hat{n}}  \,-\,
\int_{\partial{V}_i}(1-h_{i,R})\frac{\partial(1-h_{i,R})}{\partial\hat{n}} \\
&=\;\int_{\partial{V}_i}\frac{\partial{h}_{i,R}}{\partial\hat{n}}
\end{aligned}
\end{equation}
where $\hat{n}$ is the outward pointing normal of $V_i$.
Since $\partial{V}_i$ is compact and since $\partial{h}_{i,R}/\partial\hat{n}$ is uniformly bounded by the Cheng-Yau gradient estimate \cite{CY}, the Dirichlet integral $\int|\nabla{h}_{i,R}|^2$ is uniformly bounded.
Since $h_{i,R}\rightarrow{h}_i$ as $R\rightarrow\infty$ in (at least) the $C^1$ sense, we have that $\int|\nabla{h}_i|^2$ is finite by Fatou's lemma.

Let $\left<\cdot,\,\cdot\right>$ denote the $L^2$ inner product on a Riemannian manifold.
If $\eta$ is any $p$-form and $\varphi$ is a $C_c^\infty$ function, a computation gives
\be
\left|\varphi\,{d}\eta\right|^2 \,+\,\left|\varphi\,{d}^{*}\eta\right|^2
&=&\left<\varphi^2\,\eta,\,\triangle\eta\right>
\,-\,2\left<d\varphi\wedge\eta,\,\varphi\,d\eta\right> \,+\,2\left<{i}_{d\varphi}\eta,\,\varphi\,d^{*}\eta\right>.
\ee
If $\eta$ is harmonic, then by replacing $\varphi$ by $\varphi^2$ and using a H\"older inequality we easily conclude
\be
\left|\varphi\,{d}\eta\right|^2 \,+\,\left|\varphi\,{d}^{*}\eta\right|^2
&\le& 4\left|d\varphi\wedge\eta\right|^2 \,+\,4\left|{i}_{d\varphi}\eta\right|^2.
\ee
It follows that if $\eta$ is also bounded (or square-integrable, or increases like $o(r^2)$), then $d\eta=d^{*}\eta=0$ (compare with \cite{Li}, Lemma 3.1).
This is proved, in the usual way, by letting $\varphi_k$ be a cutoff function with $\varphi_k\equiv1$ in $B_p(2^k)$, $\varphi_k\equiv0$ outside $B_p(2^{k+1})$, and with $|d\varphi_k|\le2^{-k+1}$, and then sending $k\rightarrow\infty$.

Finally let $\eta=Jdh_i$.
Above we proved that $|\eta|^2=|dh_i|^2$ is integrable.
In addition, we have $\triangle{J}dh_i=J\triangle{dh_i}=0$.
This is due to the K\"ahler condition, and can be seen from the Bochner formula on 1-forms:
\be
\triangle=-\triangle_g\,+\,\Ric
\ee
(where $\triangle_g$ is the rough Laplacian), by noting that both $\triangle_g$ and $\Ric$ commute with $J$.
Therefore we have proven that $dJdh_i=d^{*}Jdh_i=0$.
\qed

\begin{lemma} \label{LemmaNonConstancy}
A function $h_i$ constructed above is non-constant provided that $(K,J,\omega_0)$ has more than one pseudoconvex boundary component, or has one pseudoconvex boundary component and at least one non-parabolic end.
\end{lemma}
\underline{\sl Pf} \; Obvious by construction. \qed

\begin{proposition} \label{PropFunctionRelations}
If the original manifold  $(K,J,\omega_0)$ has only parabolic ends or has no ends, then $0=\sum_i{J}dh_i$.
If $(K,J,\omega_0)$ has a non-parabolic end, there are no relations among the $Jdh_i$.
\end{proposition}
\underline{\sl Pf} \; Let $c=(c_1,\dots,c_l)$ and consider the function $h_c=\sum_{i=1}^lc_ih_i$ on $(K,g_0)$.
Since $h_c$ is harmonic and takes the value $c_i$ on $L_i$, it cannot be constant unless all the $c_i$ are equal.
If there is a non-parabolic end $M'$ on the original manifold, a function $\varphi'$ on $M'$ exists with $\varphi\rightarrow0$ along some subsequence of points in $M'$, and (by construction) we have $h_c<C\varphi'$ where $C=\sum{c}_i/\sup_{\partial{M'}}\varphi'$, so $h_c$ also converges to $0$ along some sequence.
If all the ends of $(K,J,\omega_0)$ are parabolic, then $h_c$ is constant when all the $c_i$ are equal.
This can be seen by noting that $(K,g_0)$ is parabolic, so $\bigcup_iL_i$ has zero capacity, which implies that the function $h_c$, being $1$ on $\bigcup_iL_i=\partial{K}$, must have zero Dirichlet integral. 
\qed

\begin{theorem} \label{ThmNonTrivDeRhamClass}
If $h_c=\sum_ic_ih_i$ is not constant and takes its maximum on $\partial{K}$, then $Jdh_c$ represents a non-trivial class in $H^1_{DR}(K)$.
\end{theorem}
\underline{\sl Pf} \; Recall that $h_i$ (if non-constant) distinguishes the boundary component $L_i$ in the sense that given any neighborhood $U_i$ of $L_i$, a number $\delta>0$ exists so that $\{h_i>1-\delta\}$ is a neighborhood of $L_i$ contained in $U_i$.
Since $h_c$ obtains its maximum on those $L_i$ for which $c_i=\sup_j\{c_j\}$, which we can take to be $1$, by Proposition \ref{PropConstructionOfFunctions} there is a number $\delta$ so that some component of $\{h_c>1-\delta\}$ is a pre-compact neighborhood that is bounded away from all other boundary components of $K$.

For convenience, give the function $h_c$ the name $x$.
For an argument by contradiction, suppose $x$ is exact, meaning a function $y:K\rightarrow\RR$ exists with $dy=-Jdx$.
Setting $z=x+\sqrt{-1}\,y$ and recalling that on functions we have $\bar\partial=\frac12\left(d+\sqrt{-1}Jd\right)$, a computation gives $\bar\partial{z}=0$, so $z:K\rightarrow\CC$ is holomorphic.
Let $V$ be a pre-compact component of $\{x\ge1-\delta\}\subset{K}$.
The image of $V$ under $z$ lies in the strip $\{1-\delta\le{x}\le1\}\subset\CC$.
By the open mapping theorem, the boundary of the image lies in the union of the lines $\{x=1-\delta\}\cup\{x=1\}$, meaning the image of $\{x\ge1-\delta\}\in\CC$ under $z$ is relatively open in $\{1-\delta\le{x}\le1\}\subset\CC$.
By continuity the image is also closed, so the image of $\{x\ge1-\delta\}\subset{K}$ is precisely $\{1-\delta\le{x}\le1\}\subset\CC$.
However, this implies that $z$ has a pole on the interior of $K$, an impossibility since both $x$ and $y$ are of class $C^1$.
\qed

\underline{\sl Proof of Theorem \ref{ThmBettiNoManifold}} \; By Theorem \ref{ThmNonTrivDeRhamClass} there is a linear map from the Hilbert space $V$ generated by $\{h_1,\dots,h_l\}$ to $H^1_{DR}(K)$.
If $V\ne\{0\}$, then Proposition \ref{PropFunctionRelations} states that the kernel is 1-dimensional if $K$ has no non-parabolic ends, and zero-dimensional if there is at least one non-parabolic end. \qed

\begin{theorem} \label{ThmBoundaryBetti}
If $h_i$ is not constant, then $b^1(L_i)\ge1$.
\end{theorem}
\underline{\sl Pf} \; Since $\nabla_0{h}_i$, the gradient of $h_i$ in the $\omega_0$-metric, is non-zero near $L_i$ by the Hopf Lemma, the isotopy lemma guarantees a $\delta_i$ so that $V_i=\{h_i\ge1-\delta_i\}$ is diffeomorphic to $L_i\times(1-\delta_i,1)$.
Thus $V_i$ and $L_i$ have the same de Rham cohomology.
The proof that $Jdh_i$ represents a nontrivial class in $H^1_{DR}(V_i)$ is identical to the proof in Theorem \ref{ThmNonTrivDeRhamClass}.
\qed

\underline{\sl Proof of Theorem \ref{ThmBettiNoEnds}} \; A condition that $h_i$ be non-constant, from Proposition \ref{PropFunctionRelations}, is that $l\ge1$ and $l+n\ge2$. \qed

\end{document}